\documentclass[12pt]{article}

\usepackage{mathrsfs}
\usepackage{amsmath}
\usepackage{amsfonts}
\usepackage{amssymb}

\usepackage{amsmath, amsthm, amsfonts, amssymb}
\def\endpf{\hbox{\vrule height1.5ex width.5em}}

\newcommand{\CC}{{\mathbb C}}
\newcommand{\RR}{{\mathbb R}}
\newcommand{\HH}{{\mathbb H}}
\newcommand{\BB}{{\mathbb B}}

\def \-{\bar}
\newcommand{\p}{\partial}

\newtheorem{theorem}{Theorem}[section]
\newtheorem{lemma}[theorem]{Lemma}
\newtheorem{corollary}[theorem]{Corollary}
\newtheorem{proposition}[theorem]{Proposition}

\newtheorem{remark}[theorem]{Remark}

\date{}

\begin{document}

\title{\bf On local holomorphic maps preserving invariant $(p, p)$-forms between bounded symmetric domains} 

\author{Yuan Yuan\footnote{ Supported in
part by National Science Foundation grant DMS-1412384}}

\vspace{3cm} \maketitle

\begin{abstract}
{\small Let $D, \Omega_1, \cdots, \Omega_m$ be irreducible bounded symmetric domains. We study local holomorphic maps from $D$ into $\Omega_1 \times\cdots
\Omega_m$ preserving the invariant $(p, p)$-forms induced from the normalized Bergman metrics up to conformal constants. We show that the local holomorphic maps extends to algebraic maps in the rank one case for any $p$ and in the rank at least two case for certain sufficiently large $p$. The total geodesy thus follows if $D=\mathbb{B}^n, \Omega_i = \mathbb{B}^{N_i}$ for any $p$ or if $D=\Omega_1 =\cdots=\Omega_m$ with rank$(D)\geq 2$ and $p$ sufficiently large. As a consequence, the algebraic correspondence between quasi-projective varieties $D / \Gamma$ preserving invariant $(p, p)$-forms is modular, where $\Gamma$ is a torsion free, discrete, finite co-volume subgroup of Aut$(D)$. This solves partially a problem raised by Mok.} 
\end{abstract}

\bigskip
\section{Introduction}

The study of local holomorphic maps preserving invariant $(p, p)$-forms was raised in \cite{Mo4}. The motivations are the geometry of local holomorphic isometric or measure-preserving maps (c.f. \cite{Mo2}\cite{Mo3}\cite{Mo4}\cite{MN}\cite{Ng1}\cite{YZ}) and the modularity of algebraic correspondence (c.f. \cite{CU} \cite{MN}).  Clozel-Ullmo \cite{CU} considered the problem of the modularity for the algebraic correspondence $Y \subset X\times X$ that commutes with a given Hecke correspondence, where $X = \Omega / \Gamma$ is a quotient of an irreducible bounded symmetric domain $\Omega$ by a torsion free, discrete subgroup $\Gamma \subset$ Aut$(\Omega)$. They reduced the problem to the characterization of algebraic correspondence preserving either the invariant metrics or the invariant volume forms. In the differential geometric formulation,  the problem can be further reduced to the local holomorphic isometries or the local holomorphic measure-preserving maps from $U \subset \Omega$ into $\Omega \times \Omega$. Clozel-Ullmo \cite{CU} proved the total geodesy for the local holomorphic isometries arised from the algebraic correspondence for $\Omega$ being the unit disc $\Delta \subset \mathbb{C}$ or a bounded symmetric domain of rank at least two, by applying the powerful metric rigidity theorem of Mok \cite{Mo1} in the latter case (see also \cite{Mo3}). Therefore, the modularity of the algebraic correspondences is derived in those cases. 

\medskip

Mok then formulates the general differential geometric problems studying the local holomorphic isometries or measure-preserving maps between bounded symmetric domains. Let $D, \Omega_1, \cdots, \Omega_m$ be irreducible bounded symmetric domains,  $U \subset D$ be a connected open set, $F: U\rightarrow \Omega_1\times \cdots \times \Omega_m$ be a holomorphic isometry up to conformal constants. Assume dim$_{\mathbb{C}} D \geq 2$. Mok proved that $F$ is totally geodesic either if $D = \Omega_1 =\cdots =\Omega_m = \mathbb{B}^n$ \cite{Mo2}; or if rank$D \geq 2$ \cite{Mo4}. Ng proved that $F$ is   totally geodesic if $m=2, D= \mathbb{B}^n, \Omega_1 = \Omega_2=\mathbb{B}^N$ with $N < 2n$ \cite{Ng1}. Yuan-Zhang proved that $F$ is totally geodesic if $D= \mathbb{B}^n, \Omega_1 =\mathbb{B}^{N_1}, \cdots, \Omega_m =\mathbb{B}^{N_m}$ without any restrictions on codimensions between unit balls in source and target domains \cite{YZ}. On the other hand, Mok constructed the non-totally geodesic embeddings from $\mathbb{B}^n$ into bounded symmetric domain $\Omega$ of rank at least two for $n$ either equal to one (see \cite{Mo3}) or greater than one \cite{Mo5}. By using the reduction of Clozel-Ullmo to holomorphic isometries, Mok deduced the modularity of algebraic correspondence between unit balls of dimension at least two \cite{Mo2}.

\medskip
  
The original modularity problem of Clozel-Ullmo by reducing to holomorphic measure-preserving maps between the irreducible bounded symmetric domain $\Omega$ and its products $\Omega \times \cdots \times\Omega$ was later solved by Mok-Ng \cite{MN}. They showed that, the holomorphic measure-preserving maps, in the case of either $\Omega = \mathbb{B}^n$ with $n \geq 2$ or rank$(\Omega) \geq 2$, are both totally geodesic.  

\medskip

The major ideas in studying the total geodesy for those local holomorphic maps are quite similar. The first step is to show that these local holomorphic maps can be extended as algebraic maps, i.e. they are defined by polynomials. Along the line of algebraicity, Mok \cite{Mo4} derived a powerful algebraicity theorem for local holomorphic isometries, that played essential roles in the proof the total geodesy. Another quite useful tool is the algebraicity of CR maps derived by Huang \cite{Hu1} (Theorem \ref{huang}), which can be used in both local holomorphic isometries  and local holomorphic measure-preserving maps .  It is now unclear to the author how the algebraicity theorem of Mok can be used for local holomorphic measure-preserving maps or local holomorphic maps preserving invariant $(p, p)$-forms. In the proof of total geodesy, if the holomorphic maps is derived from the algebraic correspondence, the Alexander type of theorems (Theorem \ref{MN} and Theorem \ref{Al}) were used in reducing the interior geometry of the holomorphic maps to the boundary geometry, combining the monodromy argument. 

\medskip

The local holomorphic maps preserving invariant $(p, p)$-forms between bounded symmetric domains were raised by Mok, as the interpolation of holomorphic isometries and holomorphic measure-preserving maps (c.f. Problem 5.3.1. in \cite{Mo4}). The total geodesy of such maps are related to the modularity of algebraic correspondence preserving invariant $(p, p)$-forms. The study of such problems follows the classical ideas by showing firstly the algebraicity and secondly the total geodesy by the monodromy argument. We invoke Huang's algebraicity theorem in the first step, while the difficulty occurs if the bounded symmetric domains are of rank at least two, i.e. the case when the holomorphic bisectional curvature has null directions. However, if $p$ is sufficiently large, more precisely, greater than the null dimensions of all irreducible components, then the  induced algebraic real hypersurfaces must have strongly pseudoconvex points. Thus Huang's algebraicity theorem can be applied. To prove the total geodesy, we use the idea developed in \cite{YZ} for the rank one case and the Alexander type theorem derived by \cite{MN} in the higher rank case. The modularity of algebraic correspondence then follows from the standard deduction in \cite{CU}.

\medskip

The geometry of local holomorphic isometries dates back to the celebrated work of Calabi \cite{Ca}. Since then, there has been quite a few works (besides the works mentioned above) on local holomorphic maps that is related to the objects considered in the current setting (e.g. \cite{DL} \cite{Eb} \cite{HY1} \cite{HY2}, etc.). The list is by no means to be complete.


  \bigskip
  
{\bf Acknowledgement}: Part of the work was done when the author
was visiting Institute of Mathematical Research at University of Hong Kong in the summer of 2012. He would like to thank Professor N. Mok for the invitation and
hospitality. The author would also like to thank Professor X. Huang, Professor N. Mok, S. Ng, Professor B. Shiffman and Yunxin Zhang for helpful discussions, and in particular, to thank Professor N. Mok for formulating Problem 5.3.1. in [Mo4] that inspires this work. 


\section{Background and main theorems}

Write $\BB^n:=\{z\in {\CC}^n: |z|<1\}$ for the unit ball in $\CC^n$.
Let $$\omega_{\mathbb{B}^n} = \sum_{j,k\le
n}\frac{1}{(1-|z|^2)^2}\big((1-|z|^2)\delta_{jk}+\-z_j z_k\big)dz_j \wedge
d\-z_k$$ be the invariant K\"ahler form associated to the normalized Bergman metric. 
Let $U\subset \BB^n$ be a connected open subset.
Consider a holomorphic map
\begin{equation}
F=(F_1,\ldots, F_m): U
\rightarrow \BB^{N_1}\times \cdots
\times \BB^{N_m}
\end{equation}
that preserves invariant $(p, p)$-forms 
up to positive conformal constant $\lambda$
in the sense that
\begin{equation}\label{preserving1}
 \lambda \omega^p_{\mathbb{B}^n}=\sum_{j=1}^{m}
F_j^*(\omega^p_{\mathbb{B}^{N_j}}).
\end{equation}
Let $p \leq n \leq N_i$ for each $1 \leq i \leq m$. 
More precisely, for each $j$, $F_j$ is a holomorphic map from
$U$ to $\BB^{N_j}$. 
%
The first main theorem is the following rigidity theorem:

\begin{theorem} \label{main1}
Suppose $n\ge 2$. Assume $F$ be the above map satisfying (\ref{preserving1}), and in addition, assume that each $F_j$ is of rank at least $p$  at some point in $U$ for $1 \leq j \leq m$. Then  each $F_j$ extends to a totally geodesic holomorphic
embedding from $({\BB}^n, \omega_{\mathbb{B}^n})$ into $({\BB}^{N_j}, \omega_{\mathbb{B}^{N_j}})$ for all $1\leq i \leq m$
and thus $\lambda = m$. 
\end{theorem}

Note that a holomorphic map from ${\BB}^n$ into ${\BB}^N$ is a
 totally geodesic embedding with respect to the normalized Bergman
 metric if and only if there are a (holomorphic) automorphism $\sigma
 \in Aut({\BB}^n)$ and an automorphism $\tau \in Aut({\BB}^N)$ such
 that $\tau\circ F\circ \sigma(z)\equiv(z,0)$.

\medskip

Let $D \subset \mathbb{C}^n$ be the Harish-Chandra realization of an irreducible bounded symmetric domain. 
Let $\omega_{D} $ be the invariant K\"ahler form associated to the normalized Bergman metric. 

Given any nonzero holomorphic tangent vector $X$, define $N_X$ to be the null space of the Hermitian bilinear form $R(X, \bar X, V, \bar W )$ on the holomorphic tangent space with respect to $X$, where $R$ is the Riemannian curvature operator of the bounded symmetric domain $D$. Define $\mathcal N = \max_{X} \dim_\mathbb{C} N_X$ to be the null dimension of $D$ \cite{Mo1}.
Note that $\mathcal N$ is invariant at each point, as $D$ is symmetric. In particular,  $\mathcal N_{D} < n$ for  the bounded symmetric domain  $D$ in $\mathbb{C}^n$ and $\mathcal N_{\mathbb{B}^n}=0$ for any $n$.

Assume rank$(D)\geq 2$. Let $U\subset D$ be a connected open subset.
Consider a holomorphic map
\begin{equation}
F=(F_1,\ldots, F_m): U
\rightarrow D \times \cdots
\times D
\end{equation}
that preserves invariant $(p, p)$-forms 
up to  the positive conformal constant $\lambda$
in the sense that
\begin{equation}\label{preserving2}
 \lambda \omega^p_{D}=\sum_{j=1}^{m}
F_j^*(\omega^p_{D}).
\end{equation}
The second main theorem is the following rigidity theorem:

\begin{theorem}\label{main2}
Let $D, F$ be as above and $F$ satisfy (\ref{preserving2}). 
Assume each $F_i$ is of full rank at some point in $D$ and $p > \mathcal N_D$. Then each $F_i$ extends to the holomorphic automorphim $F_i \in Aut(D)$ for all $1\leq i \leq m$ and thus $\lambda = m$.
\end{theorem}

\medskip

Let $X := \Omega / \Gamma$ be the quotient of an irreducible bounded symmetric domain $\Omega$ by a torsion free, discrete, finite co-volume subgroups of Aut$(\Omega)$. Let $T_Y $ be the algebraic correspondence. By definition, $Y \subset X \times X$ is  an irreducible subvariety such that the canonical projections $\pi_i|_Y$ are finite surjective morphisms for $i=1,2$. Denote the degree of $\pi_i$ by $d_i$. For any fixed point $z_0 \in X$, there exists an open set $U_0 \subset X$ of $z_0$ such that the covering map $\pi: \Omega \rightarrow X$ and $\pi_2: Y \rightarrow X$ are local biholomorphisms on $\pi^{-1}(U_0)$ and $\pi_2^{-1}(U_0)$, respectively. Moreover, $\pi_1$ is also a local biholomorphism on $\pi_2^{-1}(U_0)$.
Write $\pi_2^{-1}(U_0) = V_1 \cup \cdots \cup V_{d_2}$ as disjoint union of open sets in $Y$, and $\pi_2^{-1}(z_0)=\{y_1, \cdots, y_{d_2}\}$ with $y_{j} \in V_j$. Let $W_j = \pi_1 (V_j)$ be the open set containing $z_j = \pi_1(y_j)$. In this case, $\pi_1: V_j \rightarrow W_j, \pi_2: V_j \rightarrow U_0$ are both biholomorphisms for each $1 \leq j \leq d_2$. Denote $\pi_1|_{V_j} \circ \pi_2^{-1}$ by $\theta_j: U_0 \rightarrow W_j$ for each $j$. 

We say that the algebraic correspondence $T_Y$ locally preserves the invariant $(p, p)$-forms $\omega^p_\Omega$ on $X$ induced from the universal covering $\Omega$ if for any smooth $(n-p, n-p)$-form $\beta$ with compact support in $U_0$ such that 

$$\frac{1}{d_1} \int_X T^*_Y \beta \wedge \omega^p_\Omega = \int_X  \beta \wedge \omega^p_\Omega. $$ The detailed presentation on algebraic correspondence is referred to \cite{CU}. 

The standard deduction in \cite{CU} yields:

\begin{equation}\label{deduction}
\frac{1}{d_1} \sum_{j=1}^{d_2} \theta^{*}_j \omega^p_{\Omega} = \omega^p_{\Omega} ~ \text{on} ~U_0.\end{equation}
By lifting $\omega_\Omega$ from $X$ to its universal covering $\Omega$, the holomorphic map $\pi_1 \circ \pi^{-1}_2=(\theta_1, \cdots, \theta_{d_2}): U_0 \rightarrow \Omega \times \cdots \times \Omega$ can be considered as the germs of holomorphic map preserving the invariant $(p, p)$-forms $\omega^p_\Omega$ in the sense of (\ref{deduction}). 

\medskip

When $\Omega = \mathbb{B}^n$ for $n\geq 2$, suppose (\ref{deduction}) holds for any $p$. Then each $\theta_j \in$ Aut$(\mathbb{B}^n)$ as the straightforward application of Theorem \ref{main1}. When rank$(\Omega) \geq 2$, suppose (\ref{deduction}) holds for any $p > \mathcal N_\Omega$. Then each $\theta_j \in$ Aut$(\Omega)$ as the straightforward application of Theorem \ref{main2}. In both cases, the subvariety $Y \subset X \times X$ raised from such $\theta_j$ is the modular correspondence (see, also, the argument in \cite{CU}). Therefore, we have proved the following theorem:

\begin{theorem}\label{main3}
Let $X=\Omega / \Gamma$ and $Y \subset X \times X$ be the algebraic correspondence preserving the invariant $(p, p)$-form $\omega^p_\Omega$. Then $Y$ is modular either if $\Omega = \mathbb{B}^n$ for $n \geq 2$ and $1 \leq p \leq n$; or if rank$(\Omega) \geq 2$ for $n \geq p > \mathcal N_\Omega$. 
\end{theorem}

\medskip
Let us state the main ingredients in the proofs of Theorem \ref{main1} and Theorem \ref{main2}.
The first one is the following algebraicity theorem of Huang \cite{Hu1} for CR maps.

\begin{theorem}\label{huang}(Huang \cite{Hu1})
Let $M_1 \subset \mathbb{C}^n$ and $M_2 \subset \mathbb{C}^{N}$ be real algebraic hypersurfaces with $n > 1$ and $N \geq n$. Let $p \in M_1$ be a strongly pseudoconvex point. Suppose that $h$ is a holomorphic map from an open neighborhood $U_p$ of $p$ to $\mathbb{C}^N$ such that $h(U_p \cap M_1) \subset M_2$ and $h(p)$ is also a strongly pseudoconvex point, then $h$ is Nash algebraic.
\end{theorem}

Recall that a  function $h(z,\-{z})$ is called a Nash algebraic
function over $\CC^n$ if
 there is an irreducible polynomial $P(z,\xi,X)$ in $(z,\xi,X)\in
{\CC^n}\times  {\CC^n}\times {\CC}$ with $P(z,\-{z},h(z,\-{z}))\equiv 0$ over $\CC^n$. A holomorphic map is called Nash algebraic if each component of the map is a Nash algebraic function.

\medskip

The second main ingredient is the following normal form and the deep rigidity result proved by Huang \cite{Hu2} \cite{Hu3}. Let ${\HH}^n=\{(Z,W)\in \CC^{n-1}\times \CC: \Im W- |Z|^2>0\}$ be the Siegel upper half space. 
Assign the weight of $W$ to be 2 and that of
$Z$ to be 1. Denote by $o_{wt}(k)$ terms with weighted degree higher than $k$ and by $P(k)$ a function
of weighted degree $k$. For $p_0=(Z_0,W_0) \in \partial\mathbb{H}^n$, write $\sigma^p_0: (Z, W) \rightarrow (Z+Z_0, W+W_0+ 2i Z\cdot \overline{Z_0}) $
for the standard Heisenberg translation. The
following  Cayley transformation

\begin {equation}\label{cayley}
\rho_n(Z, W) = \bigg( \frac{2Z}{1- i W},\ \frac{1+i W}{1- iW}
\bigg)
\end{equation}
 biholomorphically maps  ${\HH}^n$ to
${\BB}^n$, and biholomorphically maps $\p\HH^n$, the Heisenberg hypersurface, to $\p\BB^n
\backslash \{(0, 1) \}$. Let $F$ be a rational proper holomorphic map from $\mathbb{H}^n$ to $\mathbb{H}^N$. By a result of Cima-Suffridge
\cite{CS}, $F$ is holomorphic in a neighborhood of $\partial\mathbb{H}^n$. The following normalization theorem is proved by Huang:

\begin{theorem}\label{moser}(Huang [Hu2-3]) For any $q \in \p{\HH^n}$, there is an element $\tau\in
Aut(\HH^{N+1})$ such that the map
$F^{**}_{q}=\tau\circ F\circ
\sigma_p^0= ((f^{**}_q)_1,\cdots,(f^{**}_q)_{n-1},\phi_q^{**},
g_q^{**})= (f_q^{**},\phi_p^{**}, g_q^{**})$ takes the following normal form:
\begin{align*}
&f^{**}_q(Z,W)=Z+ \frac{i}{2}a^{(1)}(Z)W+o_{wt}(3),\\
&\phi^{**}_q(Z,W)=\phi^{(2)}(Z)+o_{wt}(2),\\
&g^{**}_q(Z,W)=W+ o_{wt}(4)
\end{align*}
with
\begin{equation}\label{moser2}(\bar Z\cdot
a^{(1)}(Z))|Z|^2=|\phi^{(2)}(Z)|^2.
\end{equation}
Writing
$(f^{**}_q)_l =Z_j+\frac{i}{2}\sum_{k=1}^{n-1}a_{lk}Z_k W+o_{wt}(3).$
In particular, if the $(n-1)\times(n-1)$
Hermitian matrix $\left(a_{lk} \right)_{1\le l,k\le n-1} \equiv 0$ for all $q \in \p\HH^n$, then  $F$ is totally geodesic.
 \end{theorem}
 
 
Let $\Xi_j = \omega_{\HH^n} - F_j^*\omega_{\HH^{N_j}}$ and write $$\Xi =\Xi_{jk}dZ_j\otimes d\bar Z_k+\Xi_{jn}dZ_j\otimes d\bar
W+ \Xi_{nj}d W\otimes d\bar Z_j+\Xi_{nn}dW\otimes d\bar W.$$ The following proposition is proved in \cite{YZ} connecting the normal form and $\Xi$.

\begin{proposition}\label{zero}(\cite{YZ})
Assume that $F=
(f_1,\ldots, f_{n-1}, \phi, g):\HH^n\rightarrow \HH^N$
is a proper rational holomorphic map that satisfies the normalization (at the
origin) in Theorem \ref{moser}. 
Then
\begin{equation}\notag
\Xi_{jk}(0)=-2i(f_k)_{Z_j W}(0)=a_{kj}.
\end{equation}
\end{proposition}

\medskip

The third main ingredients in the proofs are the following Alexander type of theorems characterizing the automorphim of an irreducible bounded symmetric domain.

\begin{theorem}\label{MN}(Mok-Ng \cite{MN}) Let $D \subset \mathbb{C}^n$ be an irreducible bounded symmetric domain of rank $\geq 2$ in its Harish-Chandra realization. Suppose $b$ be a smooth point on $\partial D$. Let $U_b \subset \mathbb{C}^n$ be an open neighborhood of $b$ in $\mathbb{C}^n$ and $f: U_b \rightarrow \mathbb{C}^n$ be an open holomorphic embedding such that $f(U_b \cap D) \subset D$ and $f(U_b \cap \partial D ) \subset \partial D$. Then, there exists an automorphism $F: D\rightarrow D$ such that $F|_{U_b \cap D} = f|_{U_b \cap D}$.
\end{theorem}

\begin{theorem}\label{Al}(Alexander \cite{Al}) Let $\mathbb{B}^n$ be the complex unit ball of complex dimension $n\geq 2$. Let $b \in \partial \mathbb{B}^n$, $U_b$ be a connected open neighborhood of $b$ in $\mathbb{C}^n$, and $f: U_b \rightarrow \mathbb{C}^n$ be a nonconstant holomorphic map such that $f(U_b \cap \partial \mathbb{B}^n) \subset \partial \mathbb{B}^n$. Then, there exists an automorphism $F: \mathbb{B}^n \rightarrow \mathbb{B}^n$ such that $F|_{U_b \cap \mathbb{B}^n} = f|_{U_b \cap \mathbb{B}^n}$.
\end{theorem}


\section{Algebraic extension}

\subsection{Bounded symmetric domains of rank one}

\begin{lemma}\label{Griffiths}
Let the Hermitian holomorphic vector bundle $E \rightarrow X$ over a complex manifold $S$ be Griffiths negative. Then $\wedge^k E$ is also Griffiths negative.
\end{lemma}

{\noindent \bf Proof of Lemma:}
Firstly, we show that $\otimes^k E$ is Griffiths negative. Inductively, it suffices to show that $E \otimes E$ is Griffiths negative. Use $\nabla, \Theta$ to denote the connection and curvature operator of a Hermitian vector bundle (following expository in \cite{De}). It follows that $$\Theta (\nabla_{E\otimes E}) = \Theta(\nabla_E) \otimes I_E + I_E \otimes \Theta(\nabla_E)$$ (see \cite{De} formula (V-4.2') on p.258), where $I_E$ is the identity matrix with rank equal to rank$(E)$. One can easily check that $E \otimes E$ is Griffiths negative by showing that for any nonzero local holomorphic section $s$, $\bigg(\Theta(\nabla_{E\otimes E})(\frac{\partial}{\partial x^i}, \frac{\partial}{\partial \bar{x}^j}, s, \overline{s})\bigg)$ is a strictly negative-definite $m^2 \times m^2$ matrix, where $\{ x_i\}$ is the holomorphic local coordinate of $X$, $m=$dim$_{\mathbb{C}} S$.

Secondly, as $\wedge^k E$ is a subbundle of $\otimes^k E$, $\wedge^k E$ is also Griffiths negative by Proposition (6.10) in \cite{De} (p.340). In fact, there is an analogue Gauss-Codazzi equation for the vector bundle,

$$\Theta_{\wedge^k E}(u, u) = \Theta_{\otimes^k E}(u, u) - |\beta \cdot u|^2,$$
where $u \in TS \otimes \wedge^k E$ and $\beta \in \wedge^{1,0}(S, Hom(\wedge^k E, \otimes^k E / \wedge^k E))$ is the second fundament form of $\wedge^k E$ in $\otimes^k E$.
\endpf

\begin{theorem}\label{algebraicity1}
Let $F: =(F_1,\ldots, F_m): U
\rightarrow \BB^{N_1}\times \cdots
\times \BB^{N_m}$ be the holomorphic map defined on $U \subset \BB^n$ that preserves invariant $(p, p)$-forms up to the positive conformal factor
$\lambda $
in the sense that
\begin{equation}\label{preserving13}
\lambda  \omega^p_{\mathbb{B}^n}=\sum_{j=1}^{m}
F_j^*(\omega^p_{\mathbb{B}^{N_j}})~\text{on}~U.
\end{equation}
Let $p \leq n \leq N_i$ for each $1 \leq i \leq m$. Assume that each $F_j$ is of rank at least $p$ at some point in $\mathbb{B}^n$ for all $j$. 
Then $F$ is Nash algebraic.
\end{theorem}

{\noindent \bf Proof of Theorem:}
Consider $S_1 \subset \wedge^p(T U)$ and $S_2 \subset  \wedge^p(T
\BB^{N_1}) \times \cdots \times \wedge^p(T \BB^{N_m})$ as follows:

\begin{equation}
S_1 := \left\{(t, \zeta)  \in \wedge^p(T U) : \lambda \omega^p_{\mathbb{B}^n}(t)(\zeta, \bar\zeta) =1 \right\},
\end{equation}
and
\begin{equation}
\begin{split}
S_2 &:=\{ (z_1, \xi_1, \cdots, z_m, \xi_m) \in  \wedge^p(T
\BB^{N_1}) \times \cdots \times \wedge^p(T \BB^{N_m}):\\
&\ \ \ \  \omega^p_{\mathbb{B}^{N_1}}(z_1)(\xi_1, \overline{\xi_1}) + \cdots +
 \omega^p_{\mathbb{B}^{N_m}}(z_m)(\xi_m, \overline{\xi_m})=1\}.
\end{split}
\end{equation}
The defining functions $\rho_1, \rho_2$ of $S_1, S_2$ are,
respectively, as follows:
$$\rho_1=\lambda \omega^p_{\mathbb{B}^n}(t)(\zeta, \zeta)-1,$$
$$\rho_2= \omega^p_{\mathbb{B}^{N_1}}(z_1)(\xi_1, \xi_1) + \cdots + \omega^p_{\mathbb{B}^{N_m}}(z_m)(\xi_m, \xi_m) - 1.$$
Here $\{t\}, \{z_i\}$ are the canonical Euclidean coordinates on $\mathbb{C}^n, \mathbb{C}^{N_i}$ respectively.
Then one can easily check that the map $(F_1, dF_1, \cdots, F_m, dF_m)$
maps $S_1$ to $S_2$ according to the equation (\ref{preserving13}). It is obvious that $S_1, S_2$ are both real algebraic hypersurfaces by the expression of the complex hyperbolic metric of the unit ball. To finish the proof the theorem, it suffices to show that $(F_1, dF_1, \cdots, F_m, dF_m)$ maps a strongly pseudoconvex point on $S_1$ to a strongly pseudoconvex point in $S_2$.
We show the strong
pseudoconvexity of $S_2$ at $Q=(0, \xi_1, \cdots, 0, \xi_m)$ as follows and the
strong pseudoconvexity of $S_1$ follows from the same computation.

Let $\{s_{K_j}\}$ be the basis of $\wedge^p T \mathbb{B}^{N_i}$ and write $\xi_j = \xi_{K_j} s_{K_j}$, where $K_j$ is the multi-index $(k_{j_1}, \cdots, k_{j_p})$. Denote $\Delta_j  = \binom{N_j}{p}$ to be the rank of the vector bundle $\wedge^p T\mathbb{B}^{N_j}$.
By applying $\p\bar\p$ to $\rho_2$ at $Q=(0, \xi_1, \cdots, 0, \xi_m)$, we
have the following Hessian matrix
\begin{equation} \mathscr{H}=
\begin{bmatrix}
 B_1 & 0 & \cdots & 0 & 0\\
 0 & C_1 &\cdots & 0 & 0 \\
 \cdots & \cdots & \cdots & \cdots & \cdots \\
 0 & 0 & \cdots & B_m & 0 \\
 0 & 0 & \cdots & 0 & C_m \\
\end{bmatrix}
\end{equation}
where $B_j, C_j, j=1,2, \cdots, m$ are function-valued matrices with the following (in)equalities:
\begin{equation}\label{B_j1}
\begin{split}
B_j:= \bigg(\p_{z_{jk}}\p_{\bar z_{
{jl}}} \rho_2\bigg) (Q)= \bigg(- \sum_{K_j, L_j=1}^{\Delta_j} \Theta_{\wedge^p T\mathbb{B}^{N_i}}(\frac{\partial}{\partial z_{jk}}, \frac{\partial}{\partial \bar z_{jl}}, s_{K_j}
\bar{s}_{L_j})(0) \xi_{K_j} \bar\xi_{L_j}\bigg) \geq \delta |\xi_{j}|^2 I_{N_j},
\end{split}
\end{equation}
\begin{equation}
\begin{split}
C_j:=\bigg(\p_{s_{K_j}}\p_{\bar s_{L_j}}\rho_2\bigg)(Q)=p! \bigg(\delta_{K_j L_j}\bigg) \geq \delta
I_{\Delta_j}, 
\end{split}
\end{equation}
at $( 0, \xi_1, \cdots, 0, \xi_m) $ for some $\delta>0$.
Here $B_j$ is positive definite (the inequality in (\ref{B_j1}) holds) because $\wedge^p T\mathbb{B}^{N_j}$ is Griffiths negative by applying Lemma \ref{Griffiths}.
This implies that $Q \in S_2$ is a strongly pseudoconvex point.

Without loss of generality,  by
composing elements from $Aut(\BB^n)$ and $Aut(\BB^{N_1})\times \cdots \times
Aut(\BB^{N_m})$, one can assume that $F(0)=0$ and the rank of each $F_j$ at 0 is at least $p$. 
Therefore, there exists $0 \not= \zeta \in \wedge^p T_0
\BB^n$, such that $d F_j(\zeta) \not =0$ for all $j$. 
After rescaling, we assume that $(0,\zeta)\in S_1$. 
 Now the theorem follows by applying the algebracity theorem
of Huang (Theorem \ref{huang}) to the map
$(F_1, dF_1, \cdots, F_m, d F_m)$ from $S_1$ into $S_2$. $\endpf$

\subsection{Bounded symmetric domains of rank at least two}

\begin{theorem}\label{algebraicity2}
Let $D, \Omega_1, \cdots, \Omega_m$ be the Harish-Chandra realization of irreducible bounded symmetric domains in $\mathbb{C}^n, \mathbb{C}^{N_1}, \cdots, \mathbb{C}^{N_m}$ respectively. 
Let $F: =(F_1,\ldots, F_m): U \subset D
\rightarrow \Omega_{1}\times \cdots
\times \Omega_{m}$ be the holomorphic map defined on $U \subset D$ that preserves invariant $(p, p)$-forms 
in the sense that
\begin{equation}\label{preserving}
\lambda  \omega^p_D=\sum_{j=1}^{m}
F_j^*(\omega^p_{\Omega_{N_j}}),~ \text{for} ~\lambda >0, 
\end{equation}
where each $F_i$ is of rank at least $p$ at some point in $D$ and $p \leq n \leq N_i$ for all $i \in \{1, \cdots, m\}$.  Assume $p > \max\{\mathcal N_D, \mathcal N_{ \Omega_1}, \cdots, \mathcal N_{\Omega_m}\}$.
Then $F$ is Nash algebraic.
\end{theorem}

{\noindent \bf Proof of Theorem:}
The proof uses the similar argument in Theorem \ref{algebraicity1} of reducing the algebraicity of holomorphic maps to the algebraicity of CR maps.
Consider $S_1 \subset \wedge^p(T U)$ and $S_2 \subset  \wedge^p(T
\Omega_1) \times \cdots \times \wedge^p(T \Omega_m)$ as follows:

\begin{equation}
S_1 := \left\{(t, \zeta)  \in  \wedge^p(T U) : \lambda \omega^p_{D}(t)(\zeta, \bar\zeta) =1 \right\},
\end{equation}
and
\begin{equation}
\begin{split}
S_2 &:=\{ ( Z^{(1)}, \xi^{(1)}, \cdots, Z^{(m)}, \xi^{(m)}) \in \wedge^p(T
\Omega_{1}) \times \cdots \times \wedge^p(T \Omega_m):\\
&\ \ \ \  \omega^p_{\Omega_1}(Z^{(1)})(\xi^{(1)}, \overline{\xi^{(1)}}) + \cdots +
\omega^p_{\Omega_m}(Z^{(m)})(\xi^{(m)}, \overline{\xi^{(m)}}) =1\}.
\end{split}
\end{equation}
The defining functions $\rho_1, \rho_2$ of $S_1, S_2$ are,
respectively, as follows:
$$\rho_1= \lambda \omega^p_D(t)(\zeta, \overline{\zeta})-1,$$
$$\rho_2= \omega^p_{\Omega_1}(Z^{(1)})(\xi^{(1)}, \overline{\xi^{(1)}}) + \cdots + \omega^p_{\Omega_m}(Z^{(m)})(\xi^{(m)}, \overline{\xi^{(m)}})  - 1.$$
Here $\{t\}, \{Z^{(i)}\}$ are the canonical Euclidean coordinates on $\mathbb{C}^n, \mathbb{C}^{N_i}$ respectively. It is obvious that $S_1, S_2$ are both real algebraic hypersurfaces by the expression of the Bergman metrics \cite{FK}.
Moreover one can easily check that the map $( F_1, dF_1, \cdots, F_m, dF_m)$
maps $S_1$ to $S_2$ according to the equation (\ref{preserving}). To finish the proof of the theorem, it suffices to show that $(F_1, dF_1, \cdots, F_m, dF_m)$ maps a strongly pseudoconvex point on $S_1$ to a strongly pseudoconvex point in $S_2$.
We show the strong
pseudoconvexity of $S_2$ at $Q=( 0, \xi^{(1)}, \cdots, 0, \xi^{(m)})$ as follows and the
strong pseudoconvexity of $S_1$ follows from the same computation.

Let $\{s^{(j)}\}$ be the basis of $\wedge^p T \Omega_j$ and write $\xi^{(j)} = u^{(j)}_{j_1 \cdots j_p} s^{(j)}_{j_1 \cdots j_p}$, where $s^{(j)}_{j_1 \cdots j_p} = \frac{\partial}{\partial Z^{(j)}_{j_1}} \wedge \cdots \wedge  \frac{\partial}{\partial Z^{(j)}_{j_p}}$ for $j_1 < \cdots < j_p$. Note that $$\omega^p_{\Omega_j}(Z^{(j)}) (s^{(j)}_{j_1 \cdots j_p}, \overline{s^{(j)}_{j'_1 \cdots j'_p}}) = \det (\omega^{(j)}_{j_s j'_t})_{(j_1, \cdots, j_p; j'_1, \cdots, j'_p)},$$
where $\omega^{(j)}_{j_s j'_t} = \omega_{\Omega_j}(\frac{\partial}{\partial Z^{(j)}_{j_s}}, \overline{\frac{\partial}{\partial Z^{(j)}_{j'_t}}}).$ It follows by straightforward calculation under the normal coordinates that, at $0 \in \Omega_j$,

$$\partial_{Z^{(j)}_k} \bar \partial_{Z^{(j)}_l} \omega^p_{\Omega_j}(Z^{(j)}) (s^{(j)}_{j_1 \cdots j_p}, \overline{s^{(j)}_{j'_1 \cdots j'_p}}) = \begin{cases}
0 & \text{if} ~\sharp(\{j_1, \cdots, j_p\} \cap \{j'_1, \cdots, j'_p\}) \leq p-2\\
\pm R^{(j)}_{k \bar l s \bar t} & \text{if} ~\sharp(\{j_1, \cdots, j_p\} \cap \{j'_1, \cdots, j'_p\}) = p-1\\
-\sum_{1 \leq i \leq q}R^{(j)}_{k \bar l j_i \overline{j_i}}& \text{if} ~\{j_1, \cdots, j_p\} = \{j'_1, \cdots, j'_p\})
\end{cases} ,$$ where $R^{(j)}_{i \bar j s \bar t}$ is the holomorphic bisectional curvature of $\omega_{\Omega_j}$, that is semi-negative definite; $\{s, t\}$ are the two distinct indices such that $\{j_1, \cdots, j_p\} \setminus \{s\} = \{j'_1, \cdots, j'_p\} \setminus \{t\}$; and the sign in the second case depends on the positions of $s, t$ in $\{j_1, \cdots, j_p\}$ and $ \{j'_1, \cdots, j'_p\}$, respectively. 

Therefore, by applying $\p\bar\p$ to $\rho_2$ at $Q=(0, \xi^{(1)}, \cdots, 0, \xi^{(m)})$, we
have the following Hessian matrix
\begin{equation} \mathscr{H}=
\begin{bmatrix}
  B_1 & 0 & \cdots & 0 & 0\\
 0 & C_1 &\cdots & 0 & 0 \\
 \cdots & \cdots & \cdots & \cdots & \cdots \\
 0 & 0 & \cdots & B_m & 0 \\
 0 & 0 & \cdots & 0 & C_m \\
\end{bmatrix},
\end{equation}
where $B_j, C_j, D_j, j=1,2, \cdots, m$ are function-valued matrices with the following expressions:
\begin{equation}\label{B_j}
\begin{split}
B_j &:= \bigg(\p_{Z^{(j)}_{k}}\bar \p_{Z^{(j)}_{
{l}}} \rho_2\bigg) (Q) \\
&= \bigg( -\sum_{j_1< \cdots < j_p}\sum_{1 \leq i \leq q}R^{(j)}_{k \bar l j_i \overline{j_i}} |u^{(j)}_{j_1 \cdots j_p}|^2 \pm  \sum_{j_1 < \cdots < j_{p-1}, s \not= t} R^{(j)}_{k \bar l s \bar t} u^{(j)}_{j_1 \cdots s \cdots j_{p_1}} \overline{u^{(j)}_{j_1 \cdots t \cdots j_{p_1}}} \bigg)(Q) ,
\end{split}
\end{equation}
and
\begin{equation}
\begin{split}
C_j:=\bigg(\p_{u^{(j)}_{j_1 \cdots j_p} }\bar \p_{ u^{(j)}_{j'_1 \cdots j'_p} }\rho_2\bigg)(Q)=\bigg( \det (\omega^{(j)}_{j_s j'_t})_{(j_1, \cdots, j_p; j'_1, \cdots, j'_p)}\bigg)(0) , \end{split}
\end{equation}
%
where $C_j$ equals the $(\Delta_j \times \Delta_j)$-matrix of $\omega^p_{\Omega_j}$ at 0, that is positive definite, and $\Delta_j = \binom{N_j}{p} $ is the rank of the holomorphic vector bundle $\wedge^p T\Omega_j$. 

To show that each $B_j$ is positive definite for all $j$, fix $V^{(j)} = (v^{(j)}_1, \cdots, v^{(j)}_{N_j})$ to be a nonzero vector in the holomorphic tangent bundle $T\Omega_{j}$ at $0 \in \Omega_{j}$. Using a change of coordinate at $0 \in \Omega_j$ after a unitary transformation with constant entries, one can assume that $R^{(j)}_{V^{(j)} \overline{V^{(j)}} s \bar t}(0)$ is a diagonal matrix. Therefore, 

\begin{equation}
\begin{split}
V^{(j)} B_j \overline{V^{(j)}}^t &= -\sum_{j_1< \cdots < j_p}\sum_{1 \leq i \leq p}R^{(j)}_{V^{(j)} \overline{V^{(j)}} j_i \overline{j_i}}(0) |u^{(j)}_{j_1 \cdots j_p}|^2 \pm  \sum_{j_1 < \cdots < j_{p-1}, s \not= t} R^{(j)}_{V^{(j)} \overline{V^{(j)}} s \bar t}(0) u^{(j)}_{j_1 \cdots s \cdots j_{p_1}} \overline{u^{(j)}_{j_1 \cdots t \cdots j_{p_1}}} \\
&= -\sum_{j_1< \cdots < j_p}\sum_{1 \leq i \leq p}R^{(j)}_{V^{(j)} \overline{V^{(j)}} j_i \overline{j_i}}(0) |u^{(j)}_{j_1 \cdots j_p}|^2
\end{split}.
\end{equation}

Without loss of generality, by
composing with the automorphism groups, one can assume that $F(0)=0$ and also $F_1, \cdots, F_m$ are of rank at least $p$ at $0$. 
Therefore, there exists $(0, \zeta) \in S_1$, such that $d F_j(\zeta) \not =0$ for all $j$. 
Moreover, there exists at least one $(j_1, \cdots, j_p)$ such that $u^{(j)}_{j_1 \cdots j_p}(0, dF_j(\zeta)) \not= 0$. It follows that for this particular $(j_1, \cdots, j_p)$,

$$V^{(j)} B_j \overline{V^{(j)}}^t  \geq -\sum_{1 \leq i \leq p}R^{(j)}_{V^{(j)} \overline{V^{(j)}} j_i \overline{j_i}}(0) |u^{(j)}_{j_1 \cdots j_p}|^2>0, $$
where the second inequality follows from the assumption $p > \mathcal N_{\Omega_j}$, implying\\ $-\sum_{1 \leq i \leq p}R^{(j)}_{V^{(j)} \overline{V^{(j)}} j_i \overline{j_i}}(0) >0$.
This shows the positivity of $B_j$ and thus the positivity of $\mathscr{H}$. Therefore $Q \in S_2$ is a strongly pseudoconvex point.

Now the theorem follows by applying the algebracity theorem
of Huang (Theorem \ref{huang}) to the map
$(F_1, dF_1, \cdots, F_m, d F_m)$ from $S_1$ into $S_2$. $\endpf$

\medskip

The algebraicity in the case of unit balls follows also directly from Theorem \ref{algebraicity2}.

\begin{corollary}Let $F: =(F_1,\ldots, F_m): U \subset \mathbb{B}^n
\rightarrow \mathbb{B}^{N_1}\times \cdots
\times \mathbb{B}^{N_m}$ be the holomorphic map defined on $U \subset \mathbb{B}^n$ that preserves invariant $(p, p)$-forms 
in the sense that
\begin{equation}\label{preserving}
\lambda  \omega^p_{\mathbb{B}^n}=\sum_{j=1}^{m}
F_j^*(\omega^p_{\mathbb{B}^{N_j}}),~ \text{for} ~\lambda >0 ,
\end{equation}
where none of $F_i$ is a constant map and $p \leq n \leq N_i$ for all $i \in \{1, \cdots, m\}$.  
Then $F$ is Nash algebraic.
\end{corollary}

\section{Total geodesy}

\subsection{Bounded symmetric domains of rank one}

\begin{lemma}\label{difference}
Let $F: \BB^n \rightarrow \BB^N$ be a rational, proper holomorphic map. Let $\Xi= \omega_{\mathbb{B}^n} - F^* \omega_{\mathbb{B}^N}$. Then $\Xi$ is a non-negative $(1, 1)$-form in $\mathbb{B}^n$ that is real analytic on an open neighborhood of $\overline{\mathbb{B}^n}$. Moreover, $F$ is a totally geodesic embedding if and only if either $\Xi \equiv 0$ in $\mathbb{B}^n$ or $\Xi \equiv 0$ on an open piece $V$ of $\p \mathbb{B}^n$ of real dimension $2n-1$.
\end{lemma}

{\noindent \bf Proof:} It follows from the Schwarz Lemma that $\Xi$ is non-negative and $\Xi \equiv 0$ if and only if $F$ is a totally geodesic embedding. $\Xi$ is real analytic on an open neighborhood of $\overline{\mathbb{B}^n}$ is proved in Corollary 2.3 in \cite{YZ}. Moreover, $\Xi$ is closely related to the first fundamental form of the CR map between unit spheres. By  applying Proposition \ref{zero} to the second normalization $F^{**}_z= \sigma \circ F \circ \tau^0_z$ 
for each $z \in V$, it follows that $a_{kj}=0$ for all $z$. Hence $F$ is totally geodesic by Theorem \ref{moser}.
\endpf

\medskip

The following is a trivial result in calculus.

\begin{lemma}\label{rank}
Let $V\subset \CC^{N_1}$ be a connected open set and $F=(f_1, \cdots, f_{N_2}): V \rightarrow \CC^{N_2}$ be a holomorphic map. Let $\{w_i\}^{N_1}_{i=1}$ and $\{z_k\}_{k=1}^{N_2}$ be the coordinates of $\CC^{N_1}$ and $\CC^{N_2}$ respectively.  Assume $N_1 \geq p, N_2 \geq p$. Then

\begin{equation}
\begin{split}
F^*(d z_{k_1} \wedge \cdots \wedge d z_{k_p} \wedge d \overline{z_{l_1}} \wedge \cdots \wedge d \overline{z_{l_p}}) =  \sum_{i_1 < \cdots < i_p, j_1 < \cdots < j_p} & \det\bigg(\frac{\partial (f_{k_1}, \cdots, f_{k_p})}{\partial (w_{i_1}, \cdots, w_{i_p})}\bigg) \det\bigg(\overline{\frac{\partial (f_{l_1}, \cdots, f_{l_p})}{\partial (w_{j_1}, \cdots, w_{j_p})}} \bigg) \\
&d w_{i_1} \wedge \cdots \wedge d w_{i_p} \wedge d \overline{w_{j_1}} \wedge \cdots \wedge d \overline{w_{j_p}} .
\end{split}
\end{equation}
\end{lemma}

{\noindent \bf Proof of Theorem \ref{main1}:} Let $X$ be the union of the branch varieties of $F_j$ for $1 \leq j \leq m$. Since $\dim_{\CC} X \leq n-1$, for any $Q_0 \in \CC^n \setminus X$, there is a real curve $\gamma$ connecting $Q_0$ and $U$ such that any branch of $F$ is holomorphically continued along $\gamma$ to the germ of holomorphic map at $Q_0$, still denoted by $F$. Define $E= \cup_{j=1}^m \{z \in \BB^n \setminus X \big| |F_j(z)|=1\}$ and $\dim_{\RR} E \leq 2n-1$. 

At the first step, we are going to show $\dim_{\RR} E \leq 2n-2$. Suppose not. Then there is a curve $\gamma$ connecting $U$ and a point $Q_0 \in E$ such that $\dim_{\RR} O = 2n-1$, where $O \subset E$ is an open neighborhood of $Q_0$. Moreover,  assume $\{Q_0\}= \gamma \cap E.$ Without loss of generality, assume $Q_0 \in O \subset \{z \in E \big| |F_1(z)|=1\}.$ Since equation (\ref{preserving1}) holds in a small open neighborhood of $\gamma$ for any branch of $F$ by the holomorphic continuation, one has the following equation as points $Q_s$ on $\gamma$ approach $Q_0$:

$$ \lambda \omega^p_{\mathbb{B}^n}(Q_s)=\sum_{j=1}^{m}
F_j^*(\omega^p_{\mathbb{B}^{N_j}})(Q_s).$$

Denote the coordinates of $\BB^n$ and $\BB^{N_1}$ by $\{z_i\}_{i=1}^n$ and $\{w_k\}_{k=1}^{N_1}$ respectively and $F_1=(f_1, \cdots, f_{N_1})$. It follows that 

\begin{equation}\label{comparing}
\begin{split}
\lambda \omega^p_{\mathbb{B}^n}(Q_s) & \geq F_1^*(\omega^p_{\mathbb{B}^{N_1}})(Q_s) \geq  F_1^*\ \left(\frac{\sum_{k=1}^{N_1} d w_k \wedge d \overline{w_k}}{1-|w|^2}\right)^p (Q_s)\\
&= \frac{1}{(1-|F_1(Q_s)|^2)^{p}} \sum_{ k_1 < \cdots < k_p} C_{k_1, \cdots, k_p} \Theta_{k_1, \cdots, k_p}(Q_s),
\end{split}
\end{equation}
where each $C_{k_1, \cdots, k_p}$ is the nonnegative constant coefficient in front of $d w_{k_1} \wedge \cdots \wedge d w_{k_p} \wedge d \overline{w_{k_1}} \wedge \cdots \wedge d \overline{w_{k_p}}$ of $\big(\sum_{k=1}^{N_1} d w_k \wedge d \overline{w_k}\big)^p$ and 
$$\Theta_{k_1, \cdots, k_p} =  \sum_{i_1 < \cdots < i_p, j_1 < \cdots < j_p} \det\left( \frac{\partial (f_{k_1}, \cdots, f_{k_p})}{\partial (z_{i_1}, \cdots, z_{i_p})} \right) \det\left( \overline{\frac{\partial (f_{k_1}, \cdots, f_{k_p})}{\partial (z_{j_1}, \cdots, z_{j_p})}} \right) d z_{i_1} \wedge \cdots \wedge d z_{i_p} \wedge d \overline{z_{j_1}} \wedge \cdots \wedge d \overline{z_{j_p}},$$ by Lemma \ref{rank}.
Here two $(p, p)$-forms $\alpha_1, \alpha_2$ satisfy $\alpha_1 \geq \alpha_2$ if and only if $\alpha_1 - \alpha_2$ is a nonnegative $(p, p)$-form. 

By comparing the coefficients of  $d z_{i_1} \wedge \cdots \wedge d z_{i_p} \wedge d \overline{z_{i_1}} \wedge \cdots \wedge d \overline{z_{i_p}}$ in the equation (\ref{comparing}), as $Q_s \rightarrow Q_0$, $\omega^p_{\mathbb{B}^n}(Q_s)$ is bounded and hence $\frac{1}{(1-|F_1(Q_s)|^2)^{p}} \big|\frac{\partial (f_{k_1}, \cdots, f_{k_p})}{\partial (z_{i_1}, \cdots, z_{i_p})}(Q_s)\big|^2$ is bounded. If follows that for any $k_1 < \cdots < k_p, i_1 < \cdots < i_p$ and any $Q_0 \in O$, $$\frac{\partial (f_{k_1}, \cdots, f_{k_p})}{\partial (z_{i_1}, \cdots, z_{i_p})}(Q_0)=0.$$
By the uniqueness of holomorphic functions, it follows that $$\frac{\partial (f_{k_1}, \cdots, f_{k_p})}{\partial (z_{i_1}, \cdots, z_{i_p})} \equiv 0$$ on $U$. This contradicts to the assumption that $F_1$ is of rank at least $p$.

Given any point in $\partial \BB^n \setminus X$, still denoted by $Q_0 $, it follows from the previous step that there is a curve, still denoted by $\gamma$, connecting $U$ and $Q_0$, such that $\gamma \cap E = \emptyset$. Since the equation (\ref{preserving1}) holds on a small open neighborhood $O'$ of $\gamma$ and as $Q_s \in O'$ approaches $Q_0$, the coefficients of $\omega^p_{\mathbb{B}^n}(Q_s)$ go to $+\infty$, then the coefficients of $F_j^*(\omega_{\mathbb{B}^{N_j}}^p)(Q_s)$ go to $+\infty$  for a certain $j$. This implies $|F_j(Q_0)|=1$. Hence $F_j$ maps an open piece of $\partial \BB^n$ into $\partial \BB^{N_j}$. Assume that $F_j$ maps an open subset of $\partial \BB^n$ into $\partial \BB^{N_j}$ exactly for $1 \leq j \leq m_0$ after reordering $\{1, \cdots, j\}$. By the theorems of Forstneric \cite{Fo} and Cima-Suffridge \cite{CS}, each $F_j$ extends to the unique proper holomorphic map between $\BB^n$ and $\BB^{n_j}$, which is rational, for $1 \leq j \leq m_0$. Therefore, there exists an open subset of $\partial \BB^n$, still denoted by $O$ such that $F_j (O) \subset \partial\BB^{N_j}$ for $1 \leq j \leq m_0$ and $F_j (O) \subset \BB^{N_j}$ for $m_0+1 \leq j \leq m$. Rewrite the equation (\ref{preserving1}): 

\begin{equation}
\sum_{j=1}^{m_0} \big(\omega^p_{\mathbb{B}^n} - (F_j^*\omega_{\mathbb{B}^{N_j}})^p \big) + (\lambda- m_0) \omega^p_{\mathbb{B}^n} = \sum_{j=m_0+1}^m  F_j^*(\omega_{\mathbb{B}^{N_j}} ^p)
\end{equation}
Letting $\Xi_j = \omega_{\mathbb{B}^n} - F_j^*\omega_{\mathbb{B}^{N_j}}$ and applying difference formula, it follows that 

\begin{equation}\label{chang1}
\sum_{j=1}^{m_0}  \Xi_j \wedge \left( \sum_{t=0}^{p-1} \omega^t_{\mathbb{B}^n} \wedge (F_j^*\omega_{\mathbb{B}^{N_j}})^{p-1-t} \right) + ( \lambda - m_0) \omega^p_{\mathbb{B}^n} = \sum_{j=m_0+1}^m  F_j^*(\omega_{\mathbb{B}^{N_j}} ^p).
\end{equation}
For each $j$, it follows from $\Xi_j \geq 0$ that \begin{equation}\label{compare}
\Xi_j \wedge \omega_{\mathbb{B}^n}^{p-1} \leq \Xi_j \wedge \big( \sum_{t=0}^{p-1} \omega^t_{\mathbb{B}^n} \wedge (F_j^*\omega_{\mathbb{B}^{N_j}})^{p-1-t} \big) \leq p   \Xi_j \wedge \omega_{\mathbb{B}^n}^{p-1}.\end{equation}

Now we rewrite equation (\ref{chang1}) in the coordinate of the Siegel upper half space ${\HH}^n=\{(Z,W)\in \CC^{n-1}\times \CC: \Im W- |Z|^2>0\}$.
Applying the Cayley transformation $\rho_n(Z, W)$, one can
 compute 
 the normalized Bergman metric on $\HH^n$, denoted by $\omega_{\mathbb{H}^n}$, by pulling back the normalized Bergman metricon $\BB^n$, as follows:
\begin{equation}\label{metric1}
\begin{split}
\omega_{\mathbb{H}^n}&=\sum_{j,k<n}\frac{\delta_{jk}(\Im W-|Z|^2)+\bar Z_j Z_k}{(\Im
W-|Z|^2)^2}dZ_j \wedge d\bar Z_k+\frac{dW\wedge d\overline{
W}}{4(\Im W-|Z|^2)^2}\\
&\ \ \ \ +\sum_{j<n}\frac{\bar Z_jdZ_j\wedge d\overline{W}}{2i(\Im
W-|Z|^2)^2}-\sum_{j<n}\frac{Z_j dW\wedge d\bar Z_j}{2i(\Im
W-|Z|^2)^2}.
\end{split}
\end{equation}
Note that $\omega_{\mathbb{H}^n}$ is also an invariant metric under the action of the holomorphic automorphism group of $\HH^n$. 
Still denote $\rho^{-1}(Q_0)$ by $Q_0$. Without loss of generality, one may assume that $Z(Q_0)=0$ by composing the holomorphic automorphism of $\HH^n.$ One chooses $Q_s$ such that $Z(Q_s)=0$. Hence one has:

$$\omega_{\mathbb{H}^n}(Q_s) = \sum_{k<n} \frac{1}{\Im W(Q_s)} dZ_k \wedge d \bar Z_k + \frac{1}{4 (\Im W(Q_s))^2} dW \wedge d \overline{W}.$$
%
The right hand side of equation (\ref{chang1}) is bounded. However, because the blown-up rate for $\omega_{\mathbb{H}^n}^{p-1}(Q_s)(\frac{\p}{\p Z_1} \wedge \cdots \frac{\p}{\p Z_p} \wedge \frac{\p}{\p W} \wedge \frac{\p}{\p \bar Z_1} \wedge \cdots \frac{\p}{\p \bar Z_p} \wedge\frac{\p}{\p \overline{W}})$,  and thus by equation (\ref{compare}), that for  $\Xi_j \wedge \big( \sum_{t=0}^{p-1} \omega^t_{\mathbb{H}^n} \wedge (F_j^*\omega_{\mathbb{H}^{N_j}})^{p-1-t} \big)(Q_s)(\frac{\p}{\p Z_1} \wedge \cdots \frac{\p}{\p Z_p} \wedge \frac{\p}{\p W} \wedge \frac{\p}{\p \bar Z_1} \wedge \cdots \frac{\p}{\p \bar Z_p} \wedge \frac{\p}{\p \overline{W}})$ are both 
$\frac{1}{(\Im W(Q_s))^p}$, while that of $\omega_{\mathbb{H}^n}^p(Q_s)(\frac{\p}{\p Z_1} \wedge \cdots \frac{\p}{\p Z_p} \wedge \frac{\p}{\p W} \wedge \frac{\p}{\p \bar Z_1} \wedge \cdots \frac{\p}{\p \bar Z_p} \wedge \frac{\p}{\p \overline{W}})$ is $\frac{1}{(\Im W(Q_s))^{p+1}}$, which is higher. Therefore, $m_0 = \lambda$. 
It follows that 

\begin{equation}\label{chang}
\sum_{j=m_0+1}^m F_j^*(\omega_{\mathbb{B}^{N_j}} ^p) = \sum_{j=1}^{m_0}  \Xi_j \wedge \left( \sum_{t=0}^{p-1} \omega^t_{\mathbb{B}^n} \wedge (F_j^*\omega_{\mathbb{B}^{N_j}})^{p-1-t} \right) \geq \sum_{j=1}^{m_0} \Xi_j \wedge \omega^{p-1}_{\mathbb{B}^n}.
\end{equation}

We only need consider the case $p \geq 2$ (the case $p=1$ is solved in \cite{YZ}). By the similar argument as above, it follows that $\Xi_j \equiv 0$ on an open piece of $\p\BB^n$ of real dimension $2n-1$ for each $1 \leq j \leq m_0$. Hence $F_j$ is a totally geodesic embedding for each $j$ by Lemma \ref{difference}, and therefore $\Xi_j \equiv 0$ on $\BB^n$. Now the equation (\ref{chang}) reads

\begin{equation}
 \sum_{j=m_0+1}^m  F_j^*(\omega_{\mathbb{B}^{N_j}} ^p) \equiv 0 ~\text{on}~ U,
\end{equation}
implying $F_j^*(\omega_{\mathbb{B}^{N_j}} ^p) \equiv 0$ on $U$ for $m_0+1 \leq j \leq m$. It follows from Lemma \ref{rank} that for any $1 \leq k_1 < \cdots < k_p \leq N_j, 1\leq i_1 < \cdots < i_p \leq n$, 

$$\frac{\partial (f_{j k_1}, \cdots, f_{j k_p})}{\partial (z_{i_1}, \cdots, z_{i_p})} \equiv 0,$$
where $F_j = (f_{j1}, \cdots, f_{j N_j})$ for each $m_0+1\leq j \leq m$.
This contradicts to the assumption that each $F_j$ is of rank at least $p$. Hence $m_0=m$ and the theorem is proved.
\endpf


\subsection{Bounded symmetric domains of rank at least two}

{\noindent \bf Proof of Theorem \ref{main2}:} We follow the similar argument as in the proof of Theorem \ref{main1} (see also \cite{Mo2} \cite{MN}). By Theorem \ref{algebraicity2}, $F$ is algebraic. Let $X$ be the union of the branch varieties of $F_j$ and subvarieties where $F_i$ is not of full rank  for $1 \leq j \leq m$. Since $\dim_{\CC} X \leq n-1$, for any point $Q_0 \in \CC^n \setminus X$, there is a real curve $\gamma$ connecting $Q_0$ and $U$ such that any branch of $F$ is holomorphically continued along $\gamma$ to the germ of holomorphic map at $Q_0$. Fix one branch of such map, still denoted by $F$. Fix a smooth boundary point $q \in \partial D\setminus X$ and an open neighborhood $U_q$ of $q$ outside $X$, and choose a smooth curve $\gamma(t)$ such that $\gamma(0) \in U \setminus X, \gamma(1)= q$. Moreover, the equation (\ref{preserving2}) is preserved along $\gamma(t)$ unless $\gamma(t)$ or some $F_i(\gamma(t))$ touches $\partial D$. 

By the expression of the Bergman metric $\omega_D$ of $D$ (c.f. \cite{FK}), $\omega^p_D(z)$ blows up if and only if $z \rightarrow \partial D$. When $\gamma(t) \in D$, $F^*_j \omega^p_{D} (\gamma(t))$ does not blow up by equation (\ref{preserving2}), implying $F_j(\gamma(t)) \in D$. Hence $F_j (U_q \cap D) \subset D$ for all $j$. As $\gamma(t) \rightarrow q' \in \partial D \cap U_q$,  $F^*_j \omega^p_{D} (\gamma(t))$ blows up for some $j$, implying that $F_j(q') \in \partial D$. By possibly shrinking $U_q$, there exists some $j$ such that $F_j (\partial D \cap U_q) \subset \partial D$. By Theorem \ref{MN}, $F_j \in Aut(D)$. The theorem then follows by the induction. 
$\endpf$

\begin{remark}
When $p=m$, Theorem \ref{main2} is due to Mok-Ng \cite{MN}.
\end{remark}

\begin{remark}
One can also follow exactly the same argument in \cite{MN} by using Henkin-Tumanov's characterization \cite{TK} of the automorphism of $D$. For the simplicity of argument, here we invoke the Alexander type theorem by Mok-Ng.
\end{remark}

\begin{remark}
When $D = \mathbb{B}^n$, then Theorem \ref{main2} follows also by the similar argument using Alexander's Theorem (Theorem \ref{Al}), and the theorem in this case is due to Mok \cite{Mo2} and Mok-Ng \cite{MN}, when $p=1$ and $p=n$, respectively.

\end{remark}



\begin{thebibliography}{99}

\bibitem[Al]{Al} Alexander, H.: {\em Holomorphic mappings from the ball and polydisc}, Math. Ann. 209 (1974), 249-256.


\bibitem[Ca]{Ca}Calabi, E.: {\em Isometric imbedding of complex manifolds}, Ann.
of Math. (2) 58, (1953). 1--23, MR0057000, Zbl 0051.13103.

\bibitem[CS]{CS} Cima, J. A. and Suffridge, T. J.: {\em Boundary behavior of rational
proper maps}, Duke Math. J.  60  (1990),  no. 1, 135--138, MR1047119, Zbl 0694.32016.

\bibitem[CU]{CU}Clozel, L. and Ullmo, E.: {\em Correspondances modulaires et
mesures invariantes}, J. Reine Angew. Math. 558 (2003), 47--83, MR1979182, Zbl 1042.11027.



\bibitem[De]{De}Demailly, J. P.: {Complex Analytic and Differential Geometry}, online book available on the author's webpage.


\bibitem [DL]{DL} Di Scala, A. and Loi, A.: {\em K\"ahler manifolds and their relatives}, Ann. Sc. Norm. Super. Pisa Cl. Sci. (5) 9 (2010), no. 3, 495-501.

\bibitem[Eb]{Eb} Ebenfelt, P.: {\em Local Holomorphic Isometries of a Modified Projective Space into a Standard Projective Space; Rational Conformal Factors}, arxiv:1407.7476.

\bibitem[FK]{FK}Faraut, J. and Kor\'anyi, A.: {\em Function spaces and reproducing kernels on bounded symmetric domains,} J. Funct. Anal. 88 (1990), 64-89.

\bibitem[Fo]{Fo} Forstneric, F.: {\em Extending proper holomorphic mappings of
positive codimension}, Invent. Math. \textbf{95}, 31-62 (1989), MR0969413, Zbl 0633.32017.

\bibitem[Hu1]{Hu1}Huang, X.: {\em On the mapping problem for algebraic real
hypersurfaces in the complex spaces of different dimensions}, Ann.
Inst. Fourier (Grenoble) \textbf{44} (1994), no. 2, 433--463, MR1296739, Zbl 0803.32011.

\bibitem[Hu2]{Hu2} Huang, X.: {\em On a linearity problem of proper
 holomorphic maps between balls in complex spaces of
 different dimensions}, J. Differential Geom.
 \textbf{51}(1999), 13-33, MR1703603, Zbl 1042.32008.

\bibitem[Hu3]{Hu3} Huang, X.: {\em On a semi-rigidity property for holomorphic
maps}, Asian J. Math. \textbf{7}(2003), no. 4, 463--492. (A special
issue dedicated to Y. T. Siu on the ocassion of his 60th birthday.) MR2074886, Zbl 1056.32014.

\bibitem[HJ]{HJ} Huang, X. and Ji, S.: {\em On some rigidity problems in Cauchy-Riemann analysis},
Proceedings of the International Conference on Complex Geometry and
Related Fields,  89--107, AMS/IP Stud. Adv. Math., 39, Amer. Math.
Soc., Providence, RI, 2007, MR2338621, Zbl 1130.32020.

\bibitem[HJY]{HJY} Huang, X., Ji, S. and Yin, W.: {\em On the third gap for proper holomorphic maps between balls}, Math. Ann. 358 (2014), no. 1-2, 115-142.

\bibitem[HY1]{HY1} Huang, X. and Yuan, Y.: {\em Holomorphic isometry from a K\"ahler manifold into a product of complex projective manifolds}, Geom. Funct. Anal. 24 (2014), no. 3, 854-886. 

\bibitem[HY2]{HY2} Huang, X. and Yuan, Y.: {\em Submanifolds of Hermitian symmetric spaces}, Springer Proceedings in Mathematics \& Statistics in the memory of Salah Baouendi, to appear, arXiv:1410.3556.

\bibitem[Ji]{Ji} Ji, S.: {\em Algebraicity of real analytic hypersurfaces with maximal rank}, Amer. J. Math. 124 (2002), no. 6, 1083-1102. 

\bibitem[KZ]{KZ} Kim, S.-Y. and Zaitsev, D.: {\em Rigidity of CR maps between Shilov boundaries of bounded symmetric domains}, Invent. Math. 193 (2013), no. 2, 409-437.


\bibitem[Mo1]{Mo1} Mok, N.: {Metric rigidity theorems on Hermitian locally symmetric manifolds}, Series in Pure Mathematics. 6. World Scientific Publishing Co., Inc., Teaneck, NJ, 1989. xiv+278 pp.

\bibitem[Mo2]{Mo2} Mok, N.: {\em Local holomorphic isometric embeddings arising
from correspondences in the rank-1 case}, Contemporary trends in
algebraic geometry and algebraic topology (Tianjin, 2000), 155--165,
Nankai Tracts Math., 5, World Sci. Publ., River Edge, NJ, 2002, MR1945359, Zbl 1083.32019.



\bibitem[Mo3]{Mo3} Mok, N.: {\em Geometry of holomorphic isometries and related maps between bounded domains}, Geometry and analysis. No. 2, 225-270, Adv. Lect. Math. (ALM), 18, Int. Press, Somerville, MA, 2011.

\bibitem[Mo4]{Mo4} Mok, N.: {\em Extension of germs of holomorphic isometries
up to normalizing constants with respect to the Bergman metric}, J. Eur. Math. Soc. 14 (2012), no. 5, 1617-1656.


\bibitem[Mo5]{Mo5} Mok, N.: {\em Holomorphic isometries of the complex unit ball into irreducible bounded symmetric domains}, preprint, available at http://hkumath.hku.hk/~nmok/




\bibitem[MN]{MN} Mok, N. and Ng, S.: {\em Germs of measure-preserving holomorphic maps from bounded symmetric domains to their Cartesian products}, J. Reine Angew. Math.  669 (2012), 47-73.



\bibitem[Ng1]{Ng1} Ng, S.: {\em On holomorphic isometric embeddings of the unit n-ball into products of two unit m-balls}, Math. Z. 268 (2011) no. 1-2, 347-354, MR2805439, Zbl 1225.53050.

\bibitem[Ng2]{Ng2} Ng, S.: {\em Holomorphic Double Fibration and the Mapping Problems of Classical Domains}, to appear in Int. Math. Res. Not.

\bibitem[Siu]{Siu} Siu, Y.: {\em Strong rigidity of compact quotients of exceptional bounded symmetric domains}, Duke Math. J. 48(1981), 857-871.




\bibitem[TK]{TK} Tumanov, A. E. and Khenkin, G.M.: {\em Local characterization of analytic automorphisms
of classical domains} (Russian), Dokl. Akad. Nauk SSSR 267 (1982), 796-799;
English translation: Math. Notes 32 (1982), 849?852.


\bibitem[YZ]{YZ} Yuan, Y. and Zhang, Y.: {\em Rigidity for local holomorphic isometric embeddings from $B^n$ into $B^{N_1} \times ... \times B^{N_m}$ up to conformal factors},  J. Differential Geom. 90 (2012), no. 2, 329-349.
\bigskip

 \noindent Yuan Yuan, yyuan05@syr.edu, Department of Mathematics, Syracuse University, NY 13244, USA.



\end{thebibliography}
\end{document}